\documentstyle[12pt]{article}
\newcommand{\const}{\mathop{\rm const}\limits}

\newcommand{\Var}{\mathop{\rm Var}\limits}
\newcommand{\sign}{\mathop{\rm sign}\limits}
\newcommand{\diam}{\mathop{\rm diam}\limits}

\begin{document}

\begin{center}

{\bf  FACTORABLE CONTINUITY OF RANDOM FIELDS,}\\

\vspace{4mm}

{\bf  with quantitative estimation }\\

\vspace{4mm}

{\sc Ostrovsky E., Sirota L.}\\

\normalsize

\vspace{4mm}

{\it Department of Mathematics and Statistics, Bar-Ilan University,
59200, Ramat Gan, Israel.}\\
e-mail: eugostrovsky@list.ru \\

{\it Department of Mathematics and Statistics, Bar-Ilan University,
59200, Ramat Gan, Israel.}\\
e-mail: sirota3@bezeqint.net \\

\vspace{5mm}

  {\sc Abstract.}

 \vspace{3mm}

 \end{center}

   We study in this paper the sufficient conditions  for enhanced continuity of random
 fields, i.e. such that the modulus of its continuity allows the factorable representation by the product of random variable
 on the  deterministic  module of continuity.\par
  We estimate also the ordinary and  (possible) exponential moments of these random variables. \par
  We consider also  the case of random fields with heavy tails of distribution and the so-called rectangle its continuity. \par

\vspace{4mm}

 \par

 \vspace{4mm}


\vspace{3mm}

\section{ Introduction. Notations. Statement of problem.}

\vspace{3mm}

 Let $  (X = \{x\}, d), \ d = d(x,y)  $  be compact metric space relative the distance  function $ d(\cdot, \cdot),  $
 $  (\Omega, B, {\bf P}) $ be  a (sufficiently rich) probability space,  $  \xi = \xi(x) = \xi(x, \omega), \ x \in X, \  \omega \in \Omega  $
be separable  numerical valued  random field (r.f.) (process).\par
 Define as ordinary  for each function $ f: X\to R, $ not necessary to be continuous, the modulus (module) of continuity $ \Delta(f, \delta) $

$$
\Delta(f, d, \delta) = \Delta(f,\delta)  \stackrel{def}{=} \sup_{d(x,y) \le \delta} |f(x) - f(y)|, \ \delta \ge 0. \eqno(1.1)
$$
 In what follows the value $  \delta $ belongs to the closed segment

$$
\delta \in [0, \diam(d, X) ], \hspace{6mm} \diam(d, X) \stackrel{def}{=} \sup_{x,y \in X} d(x,y).
$$

 Obviously, $   \lim_{\delta \to 0+} \Delta(f, d, \delta)  = 0  $ iff the function $  f(\cdot) $ is (uniformly) continuous. \par

 \vspace{3mm}

{\bf Definition 1.1.} The r.f.  $ \xi = \xi(x) $ is said to be {\it factorable continuous} (FC), if there exists  a continuous
non-random  non-negative function $  g = g(\delta),  $ such that $  g(0) = g(0+) = 0,  $ and finite (non-negative) random variable (r.v.)
$ \tau = \tau(\omega)  $  such that

$$
\Delta(\xi(\cdot), d, \delta)  \le \tau(\omega)\cdot g(\delta), \eqno(1.2)
$$
 (factorable inequality).\\

 Such a function $   g(\delta) $ in the relation (1.2) will be named by definition as a {\it scaling function}. \par

\vspace{3mm}

 {\bf Example 1.1.} Let $ X = \overline{N} \stackrel{def}{=} (1,2, \ldots) \cup \{\infty\} $  (extended integer positive semi - axis)
 with the ordinary distance function

$$
d_N(m,n) := |1/n - 1/m|, \ n,m < \infty; \ d_N(n, \infty) = d_N(\infty, n) := 1/n, \ n < \infty;
$$

$$
 d_N(\infty,\infty) := 0. \eqno(1.3)
$$

 Let also $  \{\xi_n \}, \ n = 1,2,3, \ldots  $ be a sequence of random variables converging to zero almost everywhere:

$$
{\bf P} \left(\lim_{n \to \infty} \xi_n = 0 \right)=1.
$$

 We extend the definition of this sequence  formally as follows: $ \xi_{\infty} = 0;  $ then the random field  (random sequence)
$  \{ \xi_n \}, \ n \in  \overline{N} $ is $  d_N  $ continuous with probability one.\par

The factorable continuity  of these sequences implies the existence of the non-negative r.v. $  \eta $ and the deterministic
sequence  $ \epsilon_n  $ converging to null for which

$$
|\xi_n| \le \eta \cdot \epsilon_n.
$$

 Evidently, if the r.f. $  \xi(x) $ is factorable continuous, then it is continuous a.e. The inverse statement is also true, see
\cite{Ostrovsky4},  \cite{Ostrovsky5}, \cite{Buldygin1}, \cite{Ostrovsky1}, chapter 4,  Theorem 4.8.2, pages 219-221.  For the
random sequences $  \{ \xi_n \} $ this conclusion there is in the textbook \cite{Kantorovich1}, chapter 2, section 3. \par

 Thus, the ordinary continuity a.e.  of the separable random field is quite equivalent to the factorable continuity. \par

 Let $  \Phi = \Phi(u), \ u \in R $ be an Orlicz function, i.e. convex, even, continuous,
twice continuous differentiable in the domain $ u \ge  2, $  strictly decreasing in the
right side half real line, and such that

$$
\Phi(0) = \Phi(0+) = 0; \ \lim_{|u| \to \infty} \Phi(u) = \infty.
$$

Assume further that  these Orlicz function  $  \Phi(\cdot) $  satisfies the so-called  $ \Delta_2 $ condition:

$$
\overline{\lim}_{u \to \infty} \frac{\Phi(2u)}{\Phi(u)} < \infty.
$$
 Briefly, $  \Phi(\cdot) \in \Delta_2.  $\par

  It is proved in addition in \cite{Ostrovsky4}, more detail investigation see in \cite{Ostrovsky1}, chapter 4,  Theorem 4.8.2, pages 219-221,
that if the r.v. $ \sup_x |\xi(x)|  $ belongs to the Orlicz space $  L(\Phi)  $ builded on the our probability space
$  (\Omega, B, {\bf P}),  $  on the other words, has a light tail of distribution,
then the function $  g(\delta) $ in (1.2) may be picked  such that the r.v. $ \tau $ belongs also to
the  at the same  Orlicz space $  L(\Phi). $ \par

The generalization of this  proposition on the case when this Orlicz's function $ \Phi(\cdot)  $ does not satisfies  the  $  \Delta_2 $
condition (1.7), is  investigated in  \cite{Ostrovsky5}. \par

\vspace{4mm}

{\bf   The aim of this report is obtaining the quantitative estimation for the correspondent
random variable and deterministic  function appearing in the factorable inequality   (1.2),
up to the extent that to obtain the not significantly improve estimations.  }\par

\vspace{4mm}

 Of course, the factorable estimate (1.2) is more convenient for application than  the classical estimate  for the
some rearrangement invariant norm $   ||   \Delta(\xi(\cdot), \delta) ||;  $ see e.g. the recent publications \cite{Hairer1}.
\cite{Ostrovsky8},  \cite{Ostrovsky9}. \par

 For instance, we intend to derive the exact asymptotical behavior for the function $  g(\delta) $ as $ \delta \to 0+  $
and the sharp classical moment estimates $  |\tau|_p $ or exponential moments $  {\bf E} \exp(\lambda \tau). $
 We will denote hereafter as  usually for arbitrary r.v. $  \zeta $

$$
|\zeta|_p  \stackrel{def}{=} \left[ {\bf E} |\zeta|^p  \right]^{1/p}, \ p \ge 1.
$$

 The Orlicz's space over our probability triplet generated by the function $  \Phi(\cdot) $ will be denoted by $ L(\Phi), $ and
the classical Luxemburg's norm of the r.v. $  \zeta $ in this Orlicz's space  will be denoted by $  ||\zeta||\Phi.  $\par

  If for example $  \Phi(u) = \Phi_p(u)= |u|^p, \ p = \const \ge 1, $ then the space $  L(\Phi_p) $ coincides with Lebesgue-Riesz
space $  L_p. $ Evidently, this function $  \Phi_p(\cdot) $ satisfies the $ \Delta_2 $ condition. \par
 The detail representation of the theory of Orlicz's spaces  may be found in the  classical monographs \cite{Krasnoselsky1},
 \cite{Rao1}, \cite{Rao2}.\\

\vspace{3mm}

\section{Previous works.}

\vspace{3mm}

  {\bf A.  Metric entropy approach.} \par

\vspace{3mm}

 Let us introduce a new norm, the so-called “moment norm”, or equally Grand Lebesgue norm,
on the set of r.v. defined on our probability space by the following way: the space $ G(\psi) $ consist,
by definition, on all the centered r.v. with finite norm

$$
||\xi||G(\psi) \stackrel{def}{=} \sup_{p \ge 2} [|\xi|_p/\psi(p)]. \eqno(2.1)
$$
 Here $  \psi = \psi(p), \ 1 \le p < b, \ b = \const \in (1, \infty]  $ is some positive monotonically increasing
 continuous defined on the open interval $   (1,b), $ bounded from below numerical function. \par

   The detail investigation of these spaces, which are named  as Grand Lebesgue Spaces (GLS),
 may be found in \cite{Kozatchenko1}, \cite{Ostrovsky1}, chapters 1,2.\par

  Define the following functions

 $$
 \phi(p) := \left[\frac{p}{\psi(p)} \right]^{-1}, \hspace{5mm} \phi^*(\lambda) := \sup_p (|\lambda| p - \phi(p)), \ \lambda \in R.
 $$
The transform (non-linear) $  \phi \to \phi^*  $ is named Legendre, or Young-Fenchel transform. \par

  Let $  b = \infty; $ it is known \cite{Kozatchenko1} that the GLS space $  G \psi $ coincides with a subspace of an {\it exponential}
Orlicz's space  relative the Young-Orlicz  function

$$
N(u) = \exp(\phi^*(u)) - 1,
$$
consisting only on all the centered (mean zero) random variables.\par

 Note that if we choose as a capacity of the function $  \psi(\cdot)  $  a {\it  degenerate} function

$$
\psi_{(r)}(p) = 1, \ p = r; \ \psi_{(r)}(p) = \infty, \ p \ne r, \eqno(2.2a)
$$
where $  r = \const \in (1,\infty), $ we conclude formally

$$
||\xi||G(\psi_r)  = |\xi|_r. \eqno(2.2b)
$$
 Thus, the classical Lebesgue-Riesz spaces $  L_r $ are particular, more precisely, extremal case of GLS.\par

 Let again $ \xi = \xi(x), \ x \in X $ be separable numerical random field.  Suppose for some number $  b > 1 $
 (finite or not)

$$
\forall p \in [1,b) \ \Rightarrow   \psi_{\xi}(p) := \sup_{x \in X} |\xi(x)|_p < \infty. \eqno(2.3)
$$
 The introduced in (2.3) function is said to be {\it  natural } function for the family of random variables $ \{  \xi(x) \}, \ x \in X.  $
It is clear that

$$
\sup_{x \in X} ||\xi(x)||G\psi_{\xi} = 1.
$$
 The {\it  natural distance  } $  d = d_{\xi} =   d_{\xi}(x,y) $ generated by r.f. $ \xi(\cdot)  $ satisfying the condition (2.3)
one can introduced by the formula

$$
 d(x,y) = d_{\xi}(x,y) =   d_{\xi}(x,y)  := || \xi(x) - \xi(y) ||G\psi_{\xi}. \eqno(2.4)
$$

\vspace{3mm}

 Let us introduce for any subset $ V, \ V \subset X $ the so-called metric
{\it entropy } $ H(V, d, \epsilon) = H(V, \epsilon) $ as a natural logarithm
of a minimal quantity $ N(V,d, \epsilon) = N(V,\epsilon) = N $
of a balls $ S(V, x, \epsilon), \ x \in V: $

$$
S(V, x, \epsilon) \stackrel{def}{=} \{y, y \in V, \ d(x,y) \le \epsilon \},
$$
which cover all the set $ V: $

$$
N = \min \{M: \exists \{x_i \}, i = 1,2,…, M, \ x_i \in V, \ V
\subset \cup_{i=1}^M S(V, x_i, \epsilon ) \},
$$
and we denote also

$$
H(V,d,\epsilon) = \log N; \ S(x_0,\epsilon) \stackrel{def}{=}
 S(X, x_0, \epsilon), \ H(d, \epsilon) \stackrel{def}{=} H(X,d,\epsilon). \eqno(2.5)
$$

 It follows from Hausdorff's theorem that
$ \forall \epsilon > 0 \ \Rightarrow H(V,d,\epsilon)< \infty $ iff the
metric space $ (V, d) $ is precompact set, i.e. is the bounded set with
compact closure.\par

\vspace{3mm}

 Denote  for the function $  \psi(p) = \psi_{\xi}(p) $

$$
v(z) := \ln \psi(1/y), \ y > 0; \ v_*(w) := \inf_{z \in (0,1)} (zw + v(z)). \eqno(2.6)
$$

 It is proved in  particular in the monograph \cite{Ostrovsky1}, p. 172-176  that

$$
|| \Delta(\xi, \delta)||G\psi_{\xi} \le 9 \cdot \int_0^{\delta} \exp [ v_*(\ln 2 + H(X, d_{\xi}, \epsilon))    ] \ d \epsilon, \eqno(2.7)
$$
if of course the integral in the right-hand side of inequality (2.7) convergent. Evidently, in this case the r.f.  $ \xi(x)  $
is continuous a.e.; moreover,

$$
\lim_{\delta \to 0+} || \Delta(\xi, \delta)||G\psi_{\xi} = 0.
$$

 Analogous result  see in \cite{Pisier1}, \cite{Ledoux1}, \cite{Adler1}, \cite{Dudley1},
\cite{Egishyants1}, \cite{Fernique3}, \cite{Ostrovsky10}, \cite{Yimin1} etc.

\vspace{4mm}

  {\bf B. More modern  majorizing measure approach. }\par

\vspace{4mm}

 We note among the previous works the articles \cite{Arnold1}, \cite{Imkeller1}; see also  the preprint
 \cite{Ostrovsky7}. \par
  It was imposed in some previous articles \cite{Kwapien1}, \cite{Bednorz1} on the function $ \Phi(\cdot) $
the following $ \nabla^2 $ condition:

$$
\Phi(x) \Phi(y) \le \Phi(K(x+y)), \ \exists K = \const \in (1,\infty), \ x,y \ge 0
$$
or equally

$$
\sup_{x,y > 0} \left[ \frac{\Phi^{-1}(xy)}{\Phi^{-1}(x) + \Phi^{-1}(y)}\right] < \infty.
\eqno(2.8)
$$
{\it  We do not suppose this condition. For instance, we can consider the function of a view $ \Phi(z) = |z|^p,  $
which does not satisfy (2.1)}. \par

 Let us introduce the following constant (more exactly, functional)

  $$
C_2 = C_2(\Phi)  = \frac{\Phi^{-1}(1)}{54 K^2},\eqno(2.9)
  $$
if there exists.  Under this assumption the distance $ d = d(x_1, x_2) $ may be constructively defined by the formula:

$$
 d(x_1, x_2) = d_{\Phi}(x_1,x_2) := || \xi(x_1) -  \xi(x_2))||\Phi.
 \eqno(2.10)
$$

 We will use  further the so-called method of majorizing (minorizing) measures. Let $  m = m(\cdot)  $ be any probability measure on the set
$  X. $  Since the function $ \Phi = \Phi(z) $ is presumed to be continuous and strictly  increasing,
it follows from the relation (2.3) that $ V(d_{\Phi}) \le 1, $ where by definition

$$
V(d):= \int_X \int_X \Phi \left[ \frac{\rho(f(x_1), f(x_2))}{d(x_1,x_2)} \right] \ m(dx_1) \  m(dx_2).
\eqno(2.11)
$$

 Let us define also the following important distance function: $ w(x_1, x_2) = $
 $$
  w(x_1, x_2; V) = w(x_1, x_2; V, m ) = w(x_1, x_2; V, m,\Phi) = w(x_1, x_2; V, m,\Phi,d) \stackrel{def}{=}
 $$

$$
 6 \int_0^{d(x_1, x_2)} \left\{ \Phi^{-1} \left[ \frac{4V}{m^2(B(r,x_1))} \right] +
\Phi^{-1}  \left[ \frac{4V}{m^2(B(r,x_2))} \right] \right\} \ dr,
\eqno(2.12)
$$
 where $ m(\cdot) $ is probabilistic Borelian measure on the set $ X.$ \par
  The triangle inequality and other properties of the distance function $ w = w(x_1, x_2) $ are proved in
\cite{Kwapien1}.\par

\vspace{3mm}

{\bf Definition 2.1. } (See  \cite{Kwapien1}). The measure $ m  $ is said to be
{\it minorizing measure } relative the distance $ d = d(x_1,x_2), $ if for each values $ x_1, x_2 \in X
\ V(d) < \infty $ and moreover $ \ w(x_1,x_2; V(d)) < \infty. $\par

\vspace{3mm}

We will denote the set of all minorizing measures on the metric set $ (X,d) $  by  $ \cal{M} = \cal{M}(X).$ \par

 Evidently, if the function $ w(x_1, x_2) $ is bounded, then the minorizing measure $  m  $ is majorizing.  Inverse
proposition is not true, see  \cite{Kwapien1}, \cite{Arnold1}. \par

\vspace{3mm}

{\bf Remark 2.1.} If the measure $  m  $ is minorizing, then

$$
w(x_n, x ; V(d)) \to 0 \ \Leftrightarrow d(x_n, x) \to 0, \ n \to \infty.
$$
 Therefore, the continuity of a function relative the distance $  d_{\Phi}  $ is equivalent to
the continuity  of this function  relative the distance $  w.  $ \par

\vspace{3mm}

{\bf Remark 2.2.}  If

$$
\sup_{x_1, x_2 \in X} w(x_1, x_2; V(d)) < \infty,
$$
then the measure $ m $ is called {\it majorizing measure.} This classical definition
with theory explanation and applications basically in the investigation of local structure
of random processes and fields  belongs to
X.Fernique \cite{Fernique2},  \cite{Fernique3} and M.Talagrand
\cite{Talagrand1}, \cite{Talagrand2}.\par
 See also \cite{Bednorz1}, \cite{Dudley1}, \cite{Ledoux1},
 \cite{Ostrovsky5}, \cite{Ostrovsky6}, \cite{Ostrovsky7}. \par

 S.Kwapien and J.Rosinsky   proved  in \cite{Kwapien1}  the following inequality:

$$
{\bf E}  \Phi \left(2 \ C_2 \sup_{t \ne s} \frac{(\xi(t)-\xi(s))}{w(t,s)} \right) \le 1 + \sup_{t \ne s}
{\bf E} \Phi \left(\frac{(\xi(t) - \xi(s))}{d(t,s)}  \right).
\eqno(2.13)
$$

  As long as we choose $ d(t,s) = d_{\Phi}(t,s),  $ we have

$$
{\bf E}  \Phi \left(2 \ C_2 \sup_{t \ne s} \frac{(\xi(t)-\xi(s))}{w(t,s)} \right) \le 2.
$$

 Recall that $ \Phi = \Phi(u) $ is convex function and $ \Phi(0) = 0;  $ following

 $$
 \Phi \left(\frac{u}{2} \right) = \Phi \left(\frac{1}{2} \cdot 0  + \frac{1}{2}\cdot u \right) \le
 \frac{1}{2}  \Phi(0) + \frac{1}{2} \Phi(u) = \frac{1}{2} \Phi(u),
 $$

 We conclude on the basis of this inequality

$$
{\bf E}  \Phi \left( C_2 \sup_{t \ne s} \frac{(\xi(t)-\xi(s))}{w(t,s)} \right) \le 1,
\eqno(2.15)
$$
or equally

$$
|| \ \sup_{w(x_1, x_2) \le \delta} (\xi(x_1) - \xi(x_2)) \ || \Phi \le \delta/C_2.
\eqno(2.15a)
$$

\vspace{3mm}

 Suppose now the measure $ m $ and certain distance on the set $  X  \  d = d(x_1. x_2)  $ are such that

$$
| \ \xi(x_1) -  \xi(x_2)|_p \le d(x_1,x_2), \ p = \const \ge 1,
\eqno(2.16)
$$

$$
 m^2(B(r,x)) \ge r^{\theta}/C(\theta), \ r \in [0,1], \ \theta = \const > 0, \ C(\theta) \in (0,\infty).
\eqno(2.17)
$$

 Let also $ p = \const  > \theta. $\par

 \vspace{3mm}

  It is proved in \cite{Ostrovsky7} that for the r.f. $ \xi = \xi(x) $ the following inequality holds:
$ m \in  \cal{M} $  and

$$
|\xi(x_1) -  \xi(x_2)| \le 12 \ Z^{1/p} \ 4^{1/p} \ C^{1/p}(\theta) \  \frac{d^{1-\theta/p}(x_1, x_2)}{1-\theta/p},
\eqno(2.18)
$$
where the r.v. $  Z  $ has unit expectation:  $ {\bf E} Z = 1. $\par

 Note that in the estimate (2.18) the right-hand side is factorable. We improve its in the next section. \par

\vspace{3mm}

\section{ Main result. Ordinary Orlicz function.}

\vspace{3mm}

{\bf 1.} Suppose the continuous a.e. random field $  \xi = \xi(x), \ x \in X $ is such that for some Orlicz function grounded over source
probability space $ \Phi = \Phi(u) $ and satisfying the $  \Delta_2  $ condition

$$
|| \sup_{x \in X} \xi(x) ||\Phi < \infty. \eqno(3.1)
$$
 The sufficient conditions for the continuity of $  \xi(\cdot) $  and for the estimate (3.1) are aforementioned in the second section.\par

 Denote

$$
\theta(\delta) = \theta_{\Phi}(\delta) :=  ||\Delta(\xi, \delta) ||\Phi. \eqno(3.2)
$$
 As long as

$$
\Delta(\xi, \delta) \le 2 \sup_x |\xi(x)|, \ || \Delta(\xi, \delta)||L(\Phi) \le 2 ||\sup_x |\xi(x)| \ ||L(\Phi),
$$
 we conclude thanks to the Lebesgue dominated convergence theorem

$$
\lim_{d_{\Phi}(x,y) \to 0+} {\bf E} \Phi(|\xi(x) - \xi(y)|| = 0. \eqno(3.3)
$$
 The equality (3.3) implies the so-called on the language of the theory  of Orlicz's spaces  $  \Phi \ - $ mean
convergence. Since the Young-Orlicz function  $ \Phi(\cdot) \ $  satisfies the $  \Delta_2 \ $ condition, this
convergence is completely equivalent to the ordinary Orlicz norm space convergence, see \cite{Krasnoselsky1}, chapter 2.
Therefore

$$
\lim_{\delta \to 0+} \theta(\delta) = \lim_{\delta \to 0+} \theta_{\Phi}(\delta)  = 0. \eqno(3.4)
$$

\vspace{3mm}

{\bf 2.} We start from the relation (3.4).
Let $  a = \{  a_n \}, \ n = 1,2,\ldots  $ be arbitrary  positive strictly decreasing
numerical non - random sequence tending to zero as $  n \to \infty. $
Let also  $  b = \{  b_n \}, \ n = 1,2,\ldots  $ be arbitrary  positive strictly decreasing
numerical non-random sequence such that

$$
\sum_{n=1}^{\infty} b_n =1
$$
and such that

$$
\lim_{n \to \infty} \frac{a_n}{b_n} = 0.
$$

 Define the positive numerical values $ \delta_n  $ as a maximal solutions of the following equations

$$
\theta(\delta_n) = ||\Delta(\xi,\delta_n)||\Phi = a_n. \eqno(3.5)
$$
 Let us consider the following random variable

$$
\tau := \sum_{n=1}^{\infty} b_n \frac{\Delta(\xi, \delta_n)}{ ||\Delta(\xi, \delta_n)||\Phi  } =
 \sum_{n=1}^{\infty} b_n \frac{\Delta(\xi, \delta_n)}{ a_n }.\eqno(3.6)
$$
 The series in (3.6) converges in the $  L(\Phi) $ norm: it follows from the triangle inequality ant from the completeness
of the Orlicz spaces $  \tau \in L(\Phi) $ and

$$
||\tau||\Phi \le \sum_n b_n = 1. \eqno(3.7)
$$

 \vspace{3mm}

{\bf 3.} We deduce  from the definition (3.6) that

$$
\Delta(\xi, \delta_n) \le \tau \cdot \frac{a_n}{b_n}.
$$

 As  long as the function $  \delta \to \Delta(\xi, \delta) $ is monotonically increasing, we derive that for the arbitrary
positive value $  \delta $ there exists an unique natural value $  n = n(\delta)  $    for which $ \delta_{n+1} < \delta \le \delta_n  $
and hence

$$
\Delta(\xi, \delta) \le \tau \cdot \frac{a(n(\delta))}{b(n(\delta))}. \eqno(3.8)
$$

 The second factor in the right-hand side (3.8) tends to zero by virtue of our choosing of both the sequences
$ \{  a_n \} $  and $ \{ b_n \}. $\par

\vspace{3mm}

{\bf 4.} To obtain the {\it continuous} function  as a capacity of the scaling function
in the right-hand side of the inequality (3.8) instead $ a(n(\delta))/b(n(\delta)), $ we
introduce the following function  $ \delta \to g_1(\delta)  $  defined on some non-negative neighborhood of origin.
We define for the values $ \delta  = \delta_n, $ where $ \ n = 1,2,\ldots  $

$$
g_1(\delta) = g_1(\delta_n) \stackrel{def}{=} \frac{a(n)}{b(n)},\eqno(3.9)
$$
and define the values of this function inside the interval (open or closed) $ [\delta_{n + 1}, \ \delta_n]  $ by means of a linear
continuous interpolation (spline function). At last, put $  g_1(0) := 0; $ then
the function $ g_1(\cdot) $ is really certain scaling function, i.e. is
non-negative, strictly increasing, continuous, $ g_1(0) = g_1(0+) = 0, $ and

$$
\Delta(\xi, \delta) \le \tau \cdot g_1(\delta). \eqno(3.10)
$$

\vspace{3mm}

{\bf 5.} The function $  g_1  $ may be redefined as follows.

$$
g(\delta) :=   ||\tau||\Phi \cdot g_1(\delta),
$$
and define correspondingly the normed r.v. $ \tau_0 := \tau/||\tau||\Phi. $ It is clear that it is non-trivial:
 $  0 < ||\tau||\Phi < \infty. $ Then we have

$$
\Delta(\xi, \delta) \le \tau_0 \cdot g(\delta), \hspace{6mm}  ||\tau_0||\Phi = 1. \eqno(3.11)
$$

\vspace{3mm}

 To summarize.

 \vspace{3mm}

{\bf Theorem 3.1.} {\it Suppose that all the conditions of this section  imposed on the random field $  \xi = \xi(x) $
are satisfied. Then the modulus of its continuity $ \Delta(\xi, \delta) $ allows the factorable estimate} (3.11).\par

\vspace{4mm}

{\it  Here is an example.} \hspace{6mm} Assume that  $  \Phi(u) = \Phi_p(u) = |u|^p. $ Suppose also the continuous a.e. random field
$  \xi = \xi(x), \ x \in X $ is such that $  \sup_x |\xi(x)|_p = 1.  $ The continuity of r.f. $ \xi(x)  $ is
understanding relative the natural (finite) distance

$$
d_p(x,y) = |\xi(x) - \xi(y)|_p; \ p = \const \ge 1.
$$

{\it Denote as before }

$$
\theta_p(\delta) =  |\Delta(\xi, \delta) |_p.
$$

  The sufficient conditions for the $  d_p(\cdot, \cdot)  $ continuity of $  \xi(\cdot) $  and consistent as \\ $  \delta \to 0+ $
estimates for $  \theta_p(\delta) $   see, e.g. in \cite{Ostrovsky7}, \cite{Pisier1}. \par

 Suppose for definiteness

$$
\exists \alpha \in (0, 1] \  \Rightarrow \theta_p(\delta) \le C_1 \ \delta^{\alpha}, \ \delta \in (0,1/e).
$$

  One can choose

 $$
 b_n = \nu \ n^{-1} \ \ln^{-1 - \nu}( n + 1 ), \ \nu > 0;
 $$

$$
a_n = n^{-1 - \theta}, \ \theta > 0.
$$

  We deduce after some computations for at the same values $  \delta $

$$
\Delta(\xi, \delta) \le  C_2(\alpha,p,C_1) \cdot \tau \cdot \delta^{\alpha \ \theta/( 1 + \theta ) } \ (1 + \theta)^{-(1 + \nu)} \
|\ln \ \delta|^{1 + \nu},
$$
 where $  \tau $ in the right-hand side is unique (non-negative) random variable  for which $  {\bf E} \tau^p = 1. $ \par

We conclude further after minimisation of the right-hand side over $ \theta  $

$$
\Delta(\xi, \delta) \le  C_3(\alpha,p,C_1) \cdot \tau \cdot \delta^{\alpha }.
$$
 It is clear that the last estimate is essentially non-improvable. \par

{\bf Remark 3.1.} The optimal choosing of the sequences $  \{a_n \}, \ \{  b_n \}  $ in general case is now an open question\par

\vspace{3mm}

\section{Main result. General case of an arbitrary Orlicz function.}

\vspace{3mm}

 We do not assume this section that the Young-Orlicz function  $  \Phi $ satisfies the $  \Delta_2 $ condition.
In particular,  it may be exponential Orlicz function.\par
 We retain all the other suppositions (and notation) of previous sections. \par

 Recall in the beginning of this section then the Orlicz function $   \Psi(\cdot)  $ is called weaker than ones function
$  \Phi, $ if for all positive constant $  v; v = \const > 0 $

$$
\lim_{u \to \infty} \frac{\Psi(u v)}{\Phi(u)} = 0.
$$

Notation: $	\Psi << \Phi. $ \\

\vspace{3mm}

{\bf Theorem 4.1. } {\it We retain all the notations and conditions of the  third section.
Let also 	$ \Psi(\cdot) $ be other Orlicz function weaker than  $  \Phi(\cdot). $ Then
there exist a $ L(\Psi) \ -  $ normed r.v. $  \zeta; \ ||\zeta||\Psi = 1 $ and
the non-negative strictly increasing continuous function $ h(\cdot) $  with condition $ h(0) = h(0+) = 0, $
depending on $ \Psi, \ \Phi  $ such that}

$$
\Delta(\xi, \delta) \le \zeta  \cdot h(\delta). \eqno(4.1)
$$

\vspace{3mm}

{\bf Proof} is at the same as one for theorem 3.1. The only new feature is following. Rewrite the equality (3.3)

$$
\lim_{d_{\Phi}(x,y) \to 0+} {\bf E} \Phi(|\xi(x) - \xi(y)|) = 0. \eqno(3.3)
$$

 The equality (3.3) implies the so - called  on the language of the theory  of Orlicz's spaces  $  \Phi \ - $ mean
convergence.  But since we consider now the case when the Young-Orlicz   function  $ \Phi(\cdot) \ $ can not satisfy
 the $  \Delta_2 \ $ condition, this convergence means only that

$$
\lim_{\delta \to 0+} \theta_{\Psi}(\delta) = 0.
$$
 The scaling function $  h = h(\delta) $  may be constructed as  the function $  g = g(\delta) $ in the third section.\par
Everything else on the-still, as before. \par

\vspace{4mm}

{\bf Example 4.1.} \ Define the following {\it family } of Young-Orlicz functions

$$
\Theta_p(u) := \exp(|u|^p) - 1, \hspace{6mm} p = \const \in (0, \infty).
$$

 Note that if $ \  0 < q < p < \infty, $ then $ \Theta_q(\cdot) << \Theta_p(\cdot). $\par

  Further, let the r.f. $ \xi(x), \ x \in X  $ be a given. Suppose the r.f. $ \xi(x) $ satisfies all the conditions
of theorem 4.1 relative the Young-Orlicz function $  \Theta_p(\cdot). $ Assume also again that the number $  q  $ is arbitrary
 from the interval $   0 < q < p < \infty. $ We propose

$$
\Delta(\xi, d_{\Theta_q}, \delta) \le \tau_{\Theta_q} \cdot g_{p,q}(\delta),
$$
where the random variable $  \tau_{\Theta_q} $ belongs to the Orlicz space $  L(\Theta_q)  $ and
$  || \tau_{\Theta_q} ||\Theta_q = 1;  $ obviously, the r.v. $  \tau_{\Theta_q}  $ as well as the non-random scaling function
$ g_{p,q} = g_{p,q} (\delta) $ dependent also on the variables $  p,q.$ \par

\vspace{4mm}

{\bf Example 4.2.} \ Consider as an exception the so-called Gaussian case, i.e. when the r.f. $ \xi(x) $ is (separable)
centered Gaussian distributed. We can conclude that the correspondent Young-Orlicz function  has a form

$$
\Phi(u) = \Phi_G(u) = \exp(u^2/2) - 1.
$$

 The (centered) random variables belonging to the Orlicz space  $ L(\Phi_G)  $  are named subgaussian. \par

 The correspondent natural distance function  $  d = d(x,y) = d_G(x,y)   $ coincides  here with the $ L_2(\Omega) $ distance
between the random values $ \xi(x)  $ and $  \xi(y): $
$$
d_G(x,y) = \sqrt{ \Var(\xi(x) - \xi(y)) }.
$$

 Suppose again that the r.f. $ \xi(\cdot)  $ is $  r_G \ - $ continuous with probability one, or equally

$$
 \theta_{\Phi_G}(\delta) :=  ||\Delta(\xi, r_{G}, \delta) ||L(\Phi_G) \ \to 0, \ \delta \to 0+.
$$
 Here the $  r_G = r_G(x,y), \ x,y \in X  $ is some $ d_G $ continuous distance at the same set $  X.  $ \par

 It follows immediately from one of the main results  in the famous work of X.Fernique \cite{Fernique1} that
the module of continuity $ \Delta(\xi, r_{G}, \delta) $  allows the "good" factorization

$$
\Delta(\xi, r_{G}, \delta) \le \tau_G \cdot \tilde{g}(\delta),
$$
where $  \tilde{g}(\delta)  $ is such that $  \lim_{\delta \to 0+}\tilde{g}(\delta) = \tilde{g}(0) = 0 $
and $  \tau_G  $ is normed subgaussian: $  \tau_G \in L(\Phi_G)  $  and $  ||\tau_G|| \Phi_G = 1.  $ \par

 This means in particular

$$
{\bf P} (\tau_G > u) \le e^{ - u^2/2 }, \ u \ge 1.
$$

 On the other hands, for some positive constant $  C = C(G), \ 0 < C(G) < 1  $

$$
{\bf P} (\tau_G > u) \ge e^{ - C(G) \ u^2/2 }, \ u \ge 1.
$$

\vspace{3mm}

{\bf (Sub)-example 4.3.} \  Let in addition to the example 4.2 $ \xi(t) = w(t), \ (\xi = w), \ t \in [0, 1/e] $ be
an ordinary Brownian motion  or equally Wiener process. It is well known that here

$$
g(\delta) = \sqrt{ \delta \cdot |\ln \delta|}, \hspace{6mm} \delta \in [0, 1/e].
$$
Denote

$$
\tau_w := \sup_{ t,s \in (0,1/e) } \left[ \frac{|w(t) - w(s)|}{g(|t - s|)} \right],
$$
then

$$
(2 \pi)^{-1/2} e^{-u^2/2} \le {\bf P}(\tau_w > u) \le 4.8 \  e^{ -u^2/2 + 2 u}, \hspace{5mm} u \ge 5.
$$
see \cite{Ostrovsky11}.

\vspace{4mm}

\section{ Heavy tailed fields.}

\vspace{3mm}

  We consider in this short section the "modified" fractional continuity for the random field $ \eta = \eta(x), x \in X $  with
"very" heavy tails.  More precisely, we do not suppose that $ \forall x \in R \ \Rightarrow \eta(x) \in L_p(\Omega)  $ for some $ p \ge 1.  $
 For instance, $  \eta(\cdot) $ may be stable distributed with parameter $ \alpha, \ \alpha \in (0,2); $ cf. \cite{Marcus1},
 \cite{Nolan1}. \par
  We will  reduce the heavy tailed fields to the considered one. Namely, let us introduce  the following transformation.
   $  \xi(x) :=  Z_m(\eta(x)), $  where

$$
Z_m(y) \stackrel{def}{=} \sign(y)  \cdot [\ln( 1 + |y|)]^m, \eqno(5.1)
$$
so that $  Z_m(0) = 0, $ where the constant positive number $  m = \const  $ may be not integer.\par
 Note that the function $  Z_m(x) $ is continuous, odd, strictly increasing. \par
 Obviously, the tails of r.f. $  \xi(x) $ are much easier as ones of the r.f. $ \eta(x). $ \par
  Suppose the transformed r.f. $  \xi(x) $ satisfies all the conditions of theorems 3.1 or 4.1.  Then

 $$
 \Delta(Z_m( \eta), \delta)) \le \tau_{Z_m} \cdot g_{Z_{m}}(\delta), \eqno(5.2)
 $$
which may be interpreted as a modified (weak) factorable modulus of continuity of the heavy tailed random field $  \eta(x).  $\par

\vspace{4mm}

\section{ Concluding remarks. Rectangle continuity of random fields.}

\vspace{3mm}

 In this last section the set $  X  $ is convex  closure of open non - empty subset of whole Euclidean space  $  X = D \subset R^d. $ \par

We define as in \cite{Ral'chenko1}, \cite{Hu1} the {\it rectangle difference}  operator $ \Box[f](\vec{x}, \vec{y}) =
 \Box[f](x,y),  \ x,y \in D, \ f:D \to R $ as follows. $ \Delta^{(i)}[f](x,y) := $

$$
f(x_1,x_2, \ldots, x_{i-1}, y_i, x_{i+1}, \ldots, x_d) - f(x_1,x_2, \ldots, x_{i-1}, x_i, x_{i+1}, \ldots, x_d), \eqno(6.0)
$$
with obvious modification when $ i=1 $ or $ i=d; $
$$
\Box[f](x,y) \stackrel{def}{=} \left\{ \otimes _{i=1}^d \Delta^{(i)} \right\} [f](x,y).\eqno(6.1)
$$
 For instance, if $ d=2, $ then

$$
\Box[f](x,y) = f(y_1,y_2) - f(x_1,y_2) - f(y_1,x_2) + f(x_1,x_2).
$$

 If the function $  f: [0,1]^d \to R $ is $ d $ times continuous differentiable, then

 $$
 \Box[f](\vec{x},\vec{y}) = \int_{x_1}^{y_1}
 \int_{x_2}^{y_2} \ldots   \int_{x_d}^{y_d} \frac{ \partial^d f}{\partial x_1  \partial x_2  \ldots   \partial x_d } \ dx_1 dx_2 \ldots   dx_d .
 $$

 The {\it rectangle module of continuity  } $  \Omega(f, \vec{\delta} ) =  \Omega(f, \delta )  $  for the
(continuous  a.e.) function $ f $ and vector $ \vec{\delta} = \delta = ( \delta_1, \delta_2, \ldots, \delta_d)
\in [0,1]^d  $  may be defined as well as ordinary module of continuity $ \Delta(f,\delta) $ as follows:

$$
\Omega(f, \vec{\delta} ) \stackrel{def}{=} \sup \{ |\Box[f](x,y)|, \ (x,y): |x_i - y_i| \le \delta_i, \ i = 1,2,\ldots,d \}. \eqno(6.2)
$$

 Let $  \xi = \xi(x) = \xi(x_1, x_2, . . . , x_d) = \xi(\vec{x}),  \ \vec{x} \in D $ be a separable random field (r.f),
not necessary to be Gaussian. The sufficient condition for rectangle continuity of $ \xi(x)$ and Orlicz's norm estimates

$$
 \gamma(\xi, \delta) \stackrel{def}{=} || \Omega(\xi, \vec{\delta} )||\Phi,  \eqno(6.3)
 $$
such that

$$
\lim_{||\delta|| \to 0}\gamma(\xi, \delta) = 0  \eqno(6.4)
$$
 for it rectangle modulus of continuity are obtained in the articles \cite{Garsia1}, \cite{Ral'chenko1}, \cite{Hu1}, \cite{Ostrovsky6}. \par
 Recall that the first publication about fractional Sobolev’s inequalities \cite{Garsia1} was devoted in particular to the such a problem. \par

 It is not hard to obtain as before from (6.4) the sufficient conditions for factorable rectangle  continuity  of the r.f., i.e. the
for the estimates of a form

$$
 \Omega(\xi, \vec{\delta}) \le \nu(\omega) \cdot g(\vec{\delta}), \eqno(6.4)
$$
where $ g(\vec{\delta}) $ is continuous non-random scaling function such that

$$
\lim_{||\delta|| \to 0} g( \vec{\delta}) = 0,
$$
which may be constructed as before in the third section,
and $ \nu(\cdot) $  is random variable for which $  ||\nu||\Psi = 1. $\par


\vspace{4mm}

\begin{thebibliography}{99}

\vspace{3mm}

\bibitem{Adler1}
{\sc Adler, R. J. and Taylor, J. E.} (2007). {\it Random Fields and Geometry.} Springer, New York.

\bibitem{Arnold1}
{\sc Arnold L. and Imkeller P.} {\it On the spatial asymptotic Behavior of stochastic
Flows in Euclidean Space.} Stoch. Processes Appl., 62(1), (1996), 19-54.

\bibitem{Bednorz1}
{\sc Bednorz W. } {\it H\"older continuity of random processes.}
 arXiv:math/0703545v1 [math.PR] 19 Mar 2007.


\bibitem{Buldygin1}
{\sc Buldygin V.V.} {\it Supports of probabilistic measures in separable Banach
spaces.} Theory Probab. Appl., 1984, 29 v.3, pp. 528-532, (in Russian).

\bibitem{Dudley1}
{\sc Dudley R.M.} {\it Uniform Central Limit Theorem.}  Cambridge University Press, 1999.

\bibitem{Egishyants1}
{\sc S. A. Egishyants, E. I. Ostrovsky.}
{\it Local and global upper functions for random fields.} Teor. Veroyatnost. i Primenen., 41:4 (1996), 755-764

\bibitem{Fernique1}
 {\sc Fernique X.} (1975). {\it Regularite des trajectoires des
    function aleatiores gaussiennes.}  Ecole de Probablite de
    Saint-Flour, IV – 1974, Lecture Notes in Mathematic. {\bf 480}, 1-96, Springer Verlag, Berlin.

\bibitem{Fernique2}
{\sc Fernique X,} {\it Caracterisation de processus de trajectoires majores ou
continues.} Seminaire de Probabilit´s XII. Lecture Notes in Math. 649,
(1978), 691-706, Springer, Berlin.

\bibitem{Fernique3}
{\sc Fernique X.} {\it Regularite de fonctions aleatoires non gaussiennes.} Ecolee
de Ete de Probabilit´s de Saint-Flour XI-1981. Lecture Notes in Math.
976, (1983), 1-74, Springer, Berlin.

\bibitem{Garsia1}
{\sc A. M. Garsia, E. Rodemich, and H. Rumsey Jr.} {\it A real variable lemma and the continuity
of paths of some Gaussian processes.} Indiana Univ. Math. J. 20 (1970), no. 6, 565-578. MR0267632 (42:2534)

\bibitem{Hairer1}
{\sc Martin Hairer and Cyril Labb'e.} {\it Multiplicative stochastic heat equations on the whole space.}
arXiv:1504.07162v1 [math.AP] 27 Apr 2015

\bibitem{Hu1}
{\sc Yaozhong Hu and Khoa Le}
{\it A multiparameter Garsia-Rodemich-Rumsey inequality and some applications.}
arXiv:1211.6809v1 [math.PR] 29 Nov 2012

\bibitem{Imkeller1}
{\sc Imkeller P. and Scheutzov M.} {\it Stratonovich calculus with spatial parameters
and anticipative problem in multiplicative ergodic theory.} Ann. Probab. 27(1), (1999), 109-129.

\bibitem{Kantorovich1}
{\sc Kantorovich L.V., Akilov G.P.} {\it  Functional Analysis.}  Moskow, Nauka, (1984),  Issue 3,  (in Russian).

\bibitem{Kozatchenko1}
 {\sc Kozatchenko Yu. V., Ostrovsky E.I.} (1985). The {\it Banach Spaces of
random Variables of subgaussian type.}  Theory of Probab. and Math.  Stat.;  (in Russian). Kiev, KSU, {\bf 32}, 43-57.

\bibitem{Krasnoselsky1}
{\sc M.A.Krasnoselsky, Ya.B.Rutisky.} {\it Convex functions and Orlicz's Spaces.} P. Noordhoff LTD, The
Netherland, Groningen, 1961.

\bibitem{Kwapien1}
{\sc Kwapien S. and Rosinsky J.}  {\it Sample H¨older continuity of stochastic processes
and majorizing measures. } (2004). Seminar on Stochastic Analysis, Random
Fields and Applications IV, Progr. in Probab. {\bf 58,}  155-163. Birkh\"ouser, Basel.

\bibitem{Ledoux1}
{\sc Ledoux, M. and Talagrand, M.}  (1991). {\it Probability in Banach spaces. Isoperimetry
and processes.} Results in Math. and Rel. Areas, {\bf 3}, 23, xii+480 pp. Springer-Verlag,
Berlin. MR 1102015

\bibitem{Marcus1}
{\sc M.B Marcus, G Pisier.} {\it Some results on the continuity of stable processes and the domain of attraction of
continuous stable processes.}
Ann. Inst. H. Poincare'. Probab. Statist., 20 (1984), pp. 177-199.

\bibitem{Nolan1}
{\sc  John P Nolan. }  {\it  Continuity of symmetric stable processes.}
Journal of Multivariate Analysis,
Volume 29, Issue 1, April 1989, Pages 84-93.

\bibitem{Ostrovsky1}
 {\sc Ostrovsky E.I.} (1999). {\it Exponential estimations for Random Fields
 and its applications}, (in Russian).  Russia, Moscow-Obninsk, OINPE.

\bibitem{Ostrovsky4}
{\sc Ostrovsky E.I.} {\it About supports of  probability measures in separable Banach spaces. }
Doklady Akadeny Nauk of USSR, Band 255, (1980), $ N^o $ 6, 836-838. (in Russian).

\bibitem{Ostrovsky5}
{\sc E. Ostrovsky.} {\it Support of Borelian measures in separable Banach spaces. }
arXiv:0808.3248v1 [math.FA] 24 Aug 2008

\bibitem{Ostrovsky6}
{\sc E. Ostrovsky and L.Sirota.} {\it Continuity of functions belonging to the fractional order Sobolev-Grand Lebesgue
spaces. }
arXiv:1301.0132v1 [math.FA] 1 Jan 2013

\bibitem{Ostrovsky7}
{\sc E. Ostrovsky and L.Sirota.} {\it Simplification of the majorizing measures method, with development.}
arXiv:1302.3202v1 [math.PR] 13 Feb 2013

\bibitem{Ostrovsky8}
{\sc E. Ostrovsky and L.Sirota.} {\it Problem of Estimation of Fractional Derivative
for a Spectral Function of Gaussian Stationary Processes.}
arXiv:1504.04105v1 [math.ST] 16 Apr 2015

\bibitem{Ostrovsky9}
{\sc E. Ostrovsky and L.Sirota.} {\it Central Limit Theorem in H\"older spaces in the terms of majorizing measures.}
arXiv:1409.6054v1 [math.PR] 21 Sep 2014

\bibitem{Ostrovsky10}
{\sc E. Ostrovsky.} {\it The Anistropic Modulus of Continuity of Random Fields.}
Theory Probab. Appl., 46(1), 161-167.

\bibitem{Ostrovsky11}
{\sc E. Ostrovsky.} {\it The bilateral exponential  estimations of maximum distributions of random fields. }
 Survay of Soviet Mathematic, (1992), V.47, Issue 6 \ (288), 225-226.

\bibitem{Pisier1}
{\sc Pisier G. } {\it Conditions d'entropie assurant la continite de certains processus et  applications a' l'analyse harmonique.}
Se'minaire d'analyse fonctionelle 1979-1980, Ecole polytechnique, Paris.

\bibitem{Ral'chenko1}
{\sc Ral'chenko, K. V. } {\it The two-parameter Garsia-Rodemich-Rumsey inequality and its
application to fractional Brownian fields.} Theory Probab. Math. Statist. No. 75,
(2007), 167-178.

\bibitem{Rao1}
{\sc Rao M.M., Ren Z.D.} {\it Theory of Orlicz Spaces.} Basel-New York, Marcel
Decker, (1991).

\bibitem{Rao2}
{\sc Rao M.M., Ren Z.D.} {\it Application of Orlicz Spaces.} Basel-New York, Marcel
Decker, (2002).

\bibitem{Talagrand1}
{\sc Talagrand, M.} (1992). {\it A simple proof of the majorizing measure theorem.} Geom.
Funct. Anal. 2, no. 1, 118–125. MR 1143666

\bibitem{Talagrand2}
{\sc Talagrand, M. } (2005). {\it The generic chaining.}  Springer-Verlag,

\bibitem{Yimin1}
{\sc Yimin Xiao.} {\it Uniform Modulus of Continuity of Random Fields.} (2015),
Internet publication, Michigan State University.

\end{thebibliography}
\end{document}